\documentclass{fundam}

\publyear{25}
\papernumber{12}
\volume{195}
\issue{1-4}
\theDOI{10.46298/fi.12364}

\versionForARXIV

\usepackage{amsmath}
\usepackage{amssymb}
\usepackage[all]{xy}
\usepackage{epsfig}

\let\originalsqcap\sqcap
\let\originalsqcup\sqcup
\usepackage{mathabx}
\usepackage{stmaryrd}
\usepackage{txfonts}
\def\squareforqed{\hbox{\rlap{$\originalsqcap$}$\originalsqcup$}}

\newcommand{\kat}[1]{{\cal #1}}

\newcommand{\trans}[1]{{#1}^{\smallsmile}}

\newcommand{\lres}[2]{#1\backslash #2}
\newcommand{\rres}[2]{#1/#2}
\newcommand{\syq}[2]{{\rm syQ}(#1,#2)}

\newcommand{\id}{{\mathbb I}}

\newcommand{\epsi}{\varepsilon}
\newcommand{\dom}[1]{{\rm dom}(#1)}

\newcommand{\one}{{\i}}
\newcommand{\zero}{0}
\newcommand{\add}{{\rm add}}
\newcommand{\lbdnus}{{\rm neg}}
\newcommand{\su}{{\rm succ}}
\newcommand{\pr}{{\rm prec}}

\newcommand{\assoc}{{\rm assoc}}
\newcommand{\swap}{{\rm swap}}

\newcommand{\ubd}[2]{{\rm ubd}_{#1}(#2)}
\newcommand{\lbd}[2]{{\rm lbd}_{#1}(#2)}
\newcommand{\lub}[2]{{\rm lub}_{#1}(#2)}
\newcommand{\glb}[2]{{\rm glb}_{#1}(#2)}

\newcommand{\nepower}[1]{{\mathcal P}_{\!ne}(#1)}

\allowdisplaybreaks

\begin{document}

\title{Relation Algebraic Approach to the Real Numbers\\ The Least-Upper-Bound Property}

\author{Michael Winter\thanks{The author gratefully acknowledges
                                support from the Natural Sciences and
                                Engineering Research Council of Canada (283267).} \\
           Department of Computer Science,
           Brock University, St.\ Catharines, ON, Canada \\
           mwinter@brocku.ca
          }
          
\runninghead{M. Winter}{Relation Algebraic Approach to the Real Numbers}

\maketitle

\begin{abstract}
In this paper we continue the investigation of a real number object, i.e., an object representing the real numbers, in categories of relations. Our axiomatization is based on a relation
algebraic version of Tarski's axioms of the real numbers. It was already shown that the addition of such an object forms a dense, linear ordered abelian group. In the current paper
we will focus on the least-upper-bound property of such an object.
\end{abstract}

\section{Introduction}

In \cite{PartI} the notion of a real number object in Heyting categories with relational powers was introduced as an abstract version of the real numbers in a suitable category of relations. 
The axioms of such an object are based on Tarski's axioms of the real numbers. Due to the existence of relational powers it was possible to formulate relation-algebraic versions of Tarski's second
order axioms in a purely equational style. In addition, since Heyting algebras do not provide Boolean complements, the theory of a real number object is complement-free. All results so far
were shown without the usage of Boolean complements. Our main motivation for avoiding Boolean complements is that the results transfer immediately to so-called $L$-fuzzy relations, i.e., 
to relations that use a Heyting algebra $L$ as truth values instead of the Boolean truth values {\tt true} and {\tt false}. This makes it possible to transfer real number objects and their properties
also to the fuzzy case.

The investigation of a real number object is an important study since categories of relations are used to specify, implement, and verify programs. Usually objects of the category represent types, 
and relations represent programs of the given programming language and/or properties thereof. The current stream of investigation continued in the current paper allows to utilize the real numbers
in such an environment.

In this paper we will focus on the least-upper-bound property of a real number object. In the context of Tarski's axioms this property is essential for showing the Archimedean property of the additive group. In addition, the least-upper-bound property is also used to provide direct proofs of the Intermediate Value Theorem and/or the Heine-Borel Theorem. All these topics are planned for future work.

\section{Heyting Categories}

In this chapter we want to introduce the mathematical notions used in this paper. We start by recalling some basic notions from categories and allegories \cite{Freyd,Book,Fuzzy5}. Then 
we are going to introduce Heyting categories as an extension of division allegories defined in \cite{Freyd}, i.e., a Heyting category
is a division allegory where the lattice of relations between two objects is a Heyting algebra instead of just a distributive lattice. Heyting categories
are also a version of Dedekind categories introduced in \cite{oli80,oli95} without the requirement of completeness of lattice of relations
between two objects.

We will write $R:A\to B$ to indicate that a morphism $R$ of a category $\kat{C}$ has source $A$ and target $B$, and the collection of all
morphisms with source $A$ and target $B$ is denoted by $\kat{C}[A,B]$. 
Composition and the identity morphism are denoted by $;$ and $\id_A$, respectively. Please note that composition has
to be read from left to right, i.e., $Q;R$ means first $Q$ and then $R$.

\begin{definition}
A Heyting category $\kat{R}$ is a category satisfying the following:
\begin{enumerate}
\item For all objects $A$ and $B$ the collection $\kat{R}[A,B]$ is a Heyting algebra.
      Binary meet, binary join, relative pseudo-complement, the induced ordering, the least and the greatest element are denoted by
      $\sqcap,\sqcup,\to,\sqsubseteq,\Bot_{AB},\Top_{AB}$, respectively.
\item $Q;\Bot_{BC}=\Bot_{AC}$ for all relations $Q:A\to B$.
\item There is a monotone operation $\trans{\rule{0pt}{6pt}}$ (called converse)
      mapping a relation $Q:A\to B$ to $\trans{Q}:B\to A$
      such that for all relations $Q:A\to B$ and $R:B\to C$ the following
      holds: $\trans{(Q;R)}=\trans{R};\trans{Q}$ and $\trans{(\trans{Q})}=Q$.
\item For all relations $Q:A\to B, R:B\to C$ and $S:A\to C$ the modular inclusion
      $(Q;R)\sqcap S\sqsubseteq Q;(R\sqcap(\trans{Q};S))$ holds.
\item For all relations $R:B\to C$ and $S:A\to C$ there is a relation $S/R:A\to B$
      (called the right residual of $S$ and $R$) such that for all $X:A\to B$ the
      following holds: $X;R\sqsubseteq S\iff X\sqsubseteq \rres{S}{R}.$
\end{enumerate}
\end{definition}

Throughout the paper we will use the axioms and some basics facts such as monotonicity of the operations, $Q;\Bot_{BC}=\Bot_{AC}=\Bot_{AB};R$, $Q;(R\sqcup S)=Q;R\sqcup Q;S$,
$(Q\sqcup T);R=Q;R\sqcup T;R$, $\trans{\Bot}_{AB}=\Bot_{BA}$, and $\trans{(Q\sqcup R)}=\trans{Q}\sqcup\trans{R}$ for all $Q,T:A\to B$ and $R,S:B\to C$ without mentioning.

If we define the left residual $\lres{Q}{R}:B\to C$ of two relations $Q:A\to B$ and $R:A\to C$ by $\lres{Q}{R}=\trans{(\rres{\trans{R}}{\trans{Q}})}$
we immediately obtain $X\sqsubseteq\lres{Q}{R}$ iff $Q;X\sqsubseteq R$. Using both residuals we define the symmetric
quotient as $\syq{Q}{R}=(\lres{Q}{R})\sqcap(\rres{\trans{Q}}{\trans{R}})$. This construction is characterized
by $X\sqsubseteq\syq{Q}{R}$ iff $Q;X\sqsubseteq R$ and $R;\trans{X}\sqsubseteq Q$. 

As usual we define the pseudo-complement $Q^\star$ of a relation $Q:A\to B$ by $Q^\star=Q\to\Bot_{AB}$. Notice that the operation $^\star$ is antitone and that we have
$Q^{\star\star\star}=Q^\star$, $\Bot_{AB}^\star=\Top_{AB}$, $\Top_{AB}^\star=\Bot_{AB}$, $(Q\sqcup R)^\star=Q^\star\sqcap R^\star$ and $(Q\sqcap R)^{\star\star}=Q^{\star\star}\sqcap R^{\star\star}$.

Following usual conventions we call a relation $Q:A\to B$ complemented iff there is a relations $R:A\to B$ with $Q\sqcap R=\Bot_{AB}$ and $Q\sqcup R=\Top_{AB}$ and regular iff $Q^{\star\star}=Q$. The following lemma
relates these notions.

\begin{lemma}\label{LemCompReg}
Suppose $Q:A\to B$ has a complement $R:A\to B$. Then we have
\begin{enumerate}
	\item $Q^\star=R$ and $R^\star=Q$,
	\item $Q$ and $R$ are regular.
\end{enumerate}
\end{lemma}

\begin{shortproof}
\begin{enumerate}
	\item From $Q\sqcap R=\Bot_{AB}$ we immediately conclude $R\sqsubseteq Q^\star$. For the converse inclusion consider
		\begin{xalignat*}{2}
			R
			&= R\sqcup(Q\sqcap Q^\star)\\
			&= (R\sqcup Q)\sqcap(R\sqcup Q^\star)\\
			&= R\sqcup Q^\star, && \mbox{$R$ complement of $Q$}
		\end{xalignat*}
		which shows that $Q^\star\sqsubseteq R$. By switching the role of $Q$ and $R$ we obtain the second equation.
      \item The first property follows from $Q^{\star\star}=R^\star=Q$ using 1. The second property is shown analogously.
\end{enumerate}
\end{shortproof}

Now, we provide a version of the so-called Schr\"oder equivalences \cite{Schmidt,RedBible} known to be valid if each lattice of relations is in fact a Boolean algebra in the more general context of Heyting categories.

\begin{lemma}[Schr\"oder equivalences]\label{Lem:Schroeder}
Let $Q:A\to B$, $R:B\to C$ and $S:A\to C$ be relations. The we have
$$ Q;R\sqsubseteq S^\star \iff \trans{Q};S\sqsubseteq R^\star \iff S;\trans{R}\sqsubseteq Q^\star. $$
\end{lemma}

\begin{proof}
We only show $\Rightarrow$ of the first equivalence. All other implications follow analogously. Suppose $Q;R\sqsubseteq S^\star$. Then we have
$Q;R\sqcap S=\Bot_{AC}$. We obtain $\trans{Q};S\sqcap R\sqsubseteq \trans{Q};(S\sqcap Q;R)=\Bot_{BC}$, 
which immediately implies $\trans{Q};S\sqsubseteq R^\star$.
\end{proof}

The Schr\"oder equivalences allow to replace certain residuals by a combination of negations and compositions as the next lemma shows.

\eject

\begin{lemma}\label{Lem:ResNeg}
Suppose $Q:A\to B$, $R:A\to C$, and $S:B\to C$ are relations. Then we have
\begin{enumerate}
	\item ${\trans{Q}}^\star=\trans{{Q^\star}}$,
	\item $\rres{R^\star}{S}=(R^{\star\star};\trans{S})^\star$,
	\item $\lres{Q}{R^\star}=(\trans{Q};R^{\star\star})^\star$.
\end{enumerate}
\end{lemma}

\begin{shortproof}
\begin{enumerate}
      \item From $\trans{{Q^\star}}\sqcap\trans{Q}=\trans{(Q^\star\sqcap Q)}=\Bot_{BA}$ we obtain the inclusion $\sqsupseteq$. For the converse inclusion
      compute $\trans{{{\trans{Q}}^\star}}\sqcap Q=\trans{({\trans{Q}}^\star\sqcap\trans{Q})}=\Bot_{AB}$. This implies $\trans{{{\trans{Q}}^\star}}\sqsubseteq Q^\star$, and, hence,
      ${\trans{Q}}^\star\sqsubseteq\trans{{Q^\star}}$.
	\item We have
		\begin{xalignat*}{2}
			(\rres{R^\star}{S})\sqcap R^{\star\star};\trans{S}
			&\sqsubseteq ((\rres{R^\star}{S});S\sqcap R^{\star\star});\trans{S}\\
			&\sqsubseteq (R^\star\sqcap R^{\star\star});\trans{S}\\
			&= \Bot_{BA},
		\end{xalignat*}
		which immediately implies $\rres{R^\star}{S}\sqsubseteq(R^{\star\star};\trans{S})^\star$. For the converse inclusion we apply Lemma \ref{Lem:Schroeder}
		to the inclusion $R^{\star\star};\trans{S}\sqsubseteq(R^{\star\star};\trans{S})^{\star\star}$ and obtain $(R^{\star\star};\trans{S})^\star;S\sqsubseteq R^{\star\star\star}=R^\star$.
		This implies $(R^{\star\star};\trans{S})^\star\sqsubseteq\rres{R^\star}{S}$.
	\item We immediately compute
		\begin{xalignat*}{2}
			\lres{Q}{R^\star}
			&= \trans{(\rres{\trans{{R^\star}}}{\trans{Q}})}\\
			&= \trans{(\rres{{\trans{R}}^\star}{\trans{Q}})} && \mbox{by 1.}\\
			&= \trans{{({\trans{R}}^{\star\star};Q)^\star}} && \mbox{by 2.}\\
			&= {\trans{({\trans{R}}^{\star\star};Q)}}^\star && \mbox{by 1.}\\
			&= (\trans{Q};\trans{{{\trans{R}}^{\star\star}}})^\star\\
			&= (\trans{Q};R^{\star\star})^\star && \mbox{by 1.}\tag*{\phantom{~}}
		\end{xalignat*}
\end{enumerate}
\end{shortproof}

The following lemma summarizes some basic properties that will be used throughout the paper. A proof can be found in \cite{Schmidt,RedBible,Winter,Book}.

\begin{lemma}\label{Lem:Basics}
Let $Q:A\to B$, $R:A\to C$ and $S:C\to D$ be relations, and $i,j:A\to A$ be partial identities, i.e., $i,j\sqsubseteq\id_A$. Then we have
\begin{enumerate}
	\item $\trans{i}=i$,
	\item $i;j=i\sqcap j$ and $i;i=i$,
    	\item $(Q;\Top_{BC}\sqcap R);S=Q;\Top_{BD}\sqcap R;S$.
\end{enumerate}
\end{lemma}

The domain of a relation $R:A\to B$, i.e., the set of elements that are related to at least one other element, can be be defined as $\dom{R}=\id_A\sqcap R;\trans{R}$.
Please note that we have $\dom{R};R=R$ which we will use throughout the paper without mentioning.

A relation $Q:A\to B$ is called univalent (or partial function) iff $\trans{Q};Q\sqsubseteq\id_B$, total iff $\id_A\sqsubseteq Q;\trans{Q}$,
injective iff $\trans{Q}$ is univalent, surjective iff $\trans{Q}$ is total, a map iff $Q$ is total and univalent. Finally, a map is called a bijection (or bijective mapping) iff it is injective and surjective.
The following lemma states some
basic properties of univalent relations and maps. Again, a proof can be found in \cite{Schmidt,RedBible,Winter,Book}.

\begin{lemma}\label{Lem:Maps}
Let $f:A\to B$ be a mapping, $g:B\to A$ univalent, $Q:C\to A$, $R:C\to B$, $S:A\to C$ and $T,U:A\to D$. Then we have
\begin{enumerate}
	\item $Q;f\sqsubseteq R$ iff $Q\sqsubseteq R;\trans{f}$,
	\item $(Q;\trans{g}\sqcap R);g=Q\sqcap R;g$,
	\item $g;(T\sqcap U)=g;T\sqcap f;U$.
\end{enumerate}
\end{lemma}

The following lemma lists some important properties of the residuals and the symmetric quotient that are needed throughout the paper. A proof can be
found in \cite{Freyd,Schmidt,RedBible,Book}.

\begin{lemma}\label{Lem:ResBasics}
Suppose $Q:A\to B$, $R:A\to C$, and $S:B\to C$ are relations, and $f:D\to B$ is a map. Then we have
\begin{enumerate}
	\item $Q;(\lres{Q}{R})\sqsubseteq R$ and $(\rres{R}{S});S\sqsubseteq R$,
	\item $f;(\lres{Q}{R})=\lres{Q;\trans{f}}{R}$ and $(\rres{R}{S});\trans{f}=\rres{R}{f;S}$,
	\item $\trans{\syq{Q}{R}}=\syq{R}{Q}$,
	\item $f;\syq{Q}{R}=\syq{Q;\trans{f}}{R}$.
\end{enumerate}
\end{lemma}

For a singleton set $\{*\}$ and concrete relations we obviously have $\id_{\{*\}}=\Top_{\{*\}\{*\}}$. Furthermore, for any set $A$ the relation $\Top_{A\{*\}}$ is
actually a map. The first property together with the totality in the second property also characterize singleton sets up to isomorphism. Therefore, 
we define a unit as an abstract version of a singleton set as follows. A unit $1$ is an object so that $\id_1=\Top_{11}$ and $\Top_{A1}$ is total for every object $A$.

Considering concrete relation a map $p:1\to A$, i.e., a relation that maps $*$ to one element $a$ in $A$, can be identified with the element $a$.
Therefore we call a map $p:1\to A$ a point (of $A$).

Another important concept is the notion of a relational product, i.e., an abstract version of the Cartesian product of sets.
The object $A\times B$ is characterized by the projection relations $\pi:A\times B\to A$ and $\rho:A\times B\to B$ satisfying
$$ 
\trans{\pi};\pi\sqsubseteq\id_A,~~~
\trans{\rho};\rho\sqsubseteq\id_B,~~~
\pi;\trans{\pi}\sqcap\rho;\trans{\rho}=\id_{A\times B},~~~
\trans{\pi};\rho=\Top_{AB}. 
$$

Please note that the axioms of a relational product immediately imply that $\pi$ and $\rho$ are maps. In addition, the relational product is a product in the categorical sense in the subcategory of maps. Given relational products we will use the following abbreviations:

\begin{xalignat*}{2}
    Q\olessthan R &= Q;\trans{\pi}\sqcap R;\trans{\rho},\\
    Q\ogreaterthan S &= \pi;Q\sqcap\rho; S,\\
    Q\otimes T &= \pi;Q;\trans{\pi}\sqcap\rho;T;\trans{\rho}=Q;\trans{\pi}\ogreaterthan T;\trans{\rho}=\pi;Q\olessthan\rho;T,
\end{xalignat*}
and obtain the following properties \cite{SchmidtWinter}.

\begin{lemma}\label{Lem:Products}
If all relational products exist, then we have
\begin{enumerate}
   	\item $\trans{(Q\olessthan R)}=\trans{Q}\ogreaterthan\trans{R}$ and $\trans{(Q\ogreaterthan S)}=\trans{Q}\olessthan\trans{S}$,
  	\item If $R$ is total, then $(Q\olessthan R);\pi=Q$ and if $Q$ is total, then $(Q\olessthan R);\rho=R$,
  	\item If $S$ is surjective, then $\trans{\pi};(Q\ogreaterthan S)=Q$ and if $Q$ is surjective, then $\trans{\rho};(Q\ogreaterthan S)=S$,
	\item If $f$ is univalent, then $f;(Q\olessthan R)=f;Q\olessthan f;R$ and if $g$ is injective, then $(Q\ogreaterthan S);g=Q;g\ogreaterthan S;g$,
	\item $(Q\olessthan R);(T\ogreaterthan U)=Q;T\sqcap R;U$,
	\item $(Q\olessthan R);(T\otimes V)=Q;T\olessthan R;V$ and $(Q\otimes X);(T\ogreaterthan U)=Q;T\ogreaterthan X;U$,
	\item $Q;\trans{\pi}\ogreaterthan(R\olessthan S)=(Q\ogreaterthan R)\olessthan\rho;S$.
\end{enumerate}
\end{lemma}
We also use the following two bijective mappings $\assoc:A\times (B\times C)\to (A\times B)\times C$ and $\swap:A\times B\to B\times A$ defined by
\begin{xalignat*}{2}
	\assoc &= \pi;\trans{\pi};\trans{\pi}\sqcap\rho;\pi;\trans{\rho};\trans{\pi}\sqcap\rho;\rho;\trans{\rho}= (\id_A\otimes\pi)\olessthan\rho;\rho= \trans{\pi};\trans{\pi}\ogreaterthan(\trans{\rho}\otimes\id_C),\\
	\swap &= \pi;\trans{\rho}\sqcap\rho;\trans{\pi}= \rho\olessthan\pi= \trans{\rho}\ogreaterthan\trans{\pi}.
\end{xalignat*}

The following properties have been shown in \cite{Winter22-1}.

\begin{lemma}\label{Lem:AssocSwap}
\begin{enumerate}
	\item $\trans{\swap}=\swap$.
	\item $(Q\olessthan R);\swap=R\olessthan Q$ and $\swap;(Q\ogreaterthan S)=S\ogreaterthan Q$.
	\item $\swap;(Q\otimes T)=(T\otimes Q);\swap$.
	\item $(U\olessthan(Q\olessthan R));\assoc=(U\olessthan Q)\olessthan R$ and $\assoc;((Q\ogreaterthan S)\ogreaterthan V)=Q\ogreaterthan (S\ogreaterthan V)$.
	\item $\assoc;((Q\otimes T)\otimes X)=(Q\otimes (T\otimes X));\assoc$.
\end{enumerate}
\end{lemma}

With the maps above we are now ready to define an abelian group within a Heyting category.

\begin{definition}
A quadruple $(A,e,f,n)$ in a Heyting category $\kat{R}$ is called an abelian group iff $A$ is an object, $e:1\to A$ is a point, and $f:A\times A\to A$ and $n:A\to A$ 
are maps satisfying:
\begin{enumerate}
	\item $f$ is associative, i.e., $(\id_A\otimes f);f=\assoc;(f\otimes\id_A);f$,
	\item $e$ is the neutral element of $f$, i.e., $(\id_A\olessthan\Top_{A1};e);f=\id_A$,
	\item $n$ is the complement map, i.e., $(\id_A\olessthan n);f=\Top_{A1};e$,
	\item $f$ is commutative. i.e., $\swap;f=f$.
\end{enumerate}
\end{definition}

The next lemma lists some basic properties of abelian groups.

\begin{lemma}\label{Lem:GroupProps}
Let  $(A,e,f,n)$ be an abelian group. The we have
\begin{enumerate}
	\item $\trans{\assoc};(\id_A\otimes f);f=(f\otimes\id_A);f$,
	\item $(\Top_{A1};e\olessthan\id_A);f=\id_A$ and $(n\olessthan \id_A);f=\Top_{A1};e$,
	\item If $g,h:A\to A$ are maps with $(g\olessthan h);f=\Top_{A1};e$, then $g;n=h$,
	\item $(n\otimes n);f;n = f$.
\end{enumerate}
\end{lemma}

\begin{shortproof}
\begin{enumerate}
	\item This follows immediately from the fact that $\assoc$ is a bijective map.
	\item Both equations follow immediately form the commutativity of $f$.
	\item We immediately compute
		\begin{xalignat*}{2}
			g;n
			&= g;n;(\id_A\olessthan\Top_{A1};e);f && \mbox{$e$ is neutral}\\
			&= (g;n\olessthan\Top_{A1};e);f && \mbox{Lemma \ref{Lem:Basics}(3)}\\
			&=  (g;n\olessthan(g\olessthan h);f);f && \mbox{assumption}\\
			&=  (g;n\olessthan(g\olessthan h));(\id_A\otimes f);f && \mbox{Lemma \ref{Lem:Products}(6)}\\
			&=  (g;n\olessthan(g\olessthan h));\assoc;(f\otimes\id_A);f && \mbox{$f$ associative}\\
			&=  ((g;n\olessthan g)\olessthan h);(f\otimes\id_A);f && \mbox{Lemma \ref{Lem:AssocSwap}(4)}\\
			&=  (g;(n\olessthan\id_A)\olessthan h);(f\otimes\id_A);f && \mbox{Lemma \ref{Lem:Maps}(3)}\\
			&=  (g;(n\olessthan\id_A);f\olessthan h);f && \mbox{Lemma \ref{Lem:Products}(6)}\\
			&=  (g;\Top_{A1};e\olessthan h);f && \mbox{by 2.}\\
			&=  (\Top_{A1};e\olessthan h);f && \mbox{$g$ total}\\
			&=  h;(\Top_{A1};e\olessthan\id_A);f && \mbox{Lemma \ref{Lem:Basics}(3)}\\			
			&=  h. && \mbox{by 2.}
		\end{xalignat*}
 	\item First of all, the relation $b=(\swap\otimes\id_{A\times A});\assoc;(\trans{\assoc};\swap\otimes\id_A);\trans{\assoc}:(A\times A)\times (A\times A)\to (A\times A)\times (A\times A)$ is a bijective map because $\swap$ and $\assoc$ are.
 		Furthermore, we have
 		\begin{xalignat*}{2}
 			\lefteqn{((n\otimes n)\olessthan\id_{A\times A});(\swap\otimes\id_{A\times A});\assoc}\\
 			&= ((\pi;n\olessthan\rho;n)\olessthan\id_{A\times A});(\swap\otimes\id_{A\times A});\assoc\\
 			&= ((\pi;n\olessthan\rho;n);\swap\olessthan\id_{A\times A});\assoc && \mbox{Lemma \ref{Lem:Products}(6)}\\ 
 			&= ((\rho;n\olessthan\pi;n)\olessthan\id_{A\times A});\assoc && \mbox{Lemma \ref{Lem:AssocSwap}(2)}\\
 			&= ((\rho;n\olessthan\pi;n)\olessthan(\pi\olessthan\rho));\assoc\\
 			&= ((\rho;n\olessthan\pi;n)\olessthan\pi)\olessthan\rho && \mbox{Lemma \ref{Lem:AssocSwap}(4)}\\
 			\lefteqn{(((\rho;n\olessthan\pi;n)\olessthan\pi)\olessthan\rho);(\trans{\assoc};\swap\otimes\id_A);\trans{\assoc}}\\
 			&= (((\rho;n\olessthan\pi;n)\olessthan\pi);\trans{\assoc};\swap\olessthan\rho);\trans{\assoc} && \mbox{Lemma \ref{Lem:Products}(6)}\\
 			&= ((\rho;n\olessthan(\pi;n\olessthan\pi));\swap\olessthan\rho);\trans{\assoc} && \mbox{Lemma \ref{Lem:AssocSwap}(4)}\\
 			&= (((\pi;n\olessthan\pi)\olessthan\rho;n)\olessthan\rho);\trans{\assoc} && \mbox{Lemma \ref{Lem:AssocSwap}(2)}\\
 			&= (\pi;n\olessthan\pi)\olessthan(\rho;n\olessthan\rho) && \mbox{Lemma \ref{Lem:AssocSwap}(4)}\\
 			&= \pi;(n\olessthan\id_A)\olessthan\rho;(n\olessthan\id_A) && \mbox{Lemma \ref{Lem:Maps}(3)}\\
 			&= (n\olessthan\id_A)\otimes(n\olessthan\id_A),
 		\end{xalignat*}
 		i,e., $((n\otimes n)\olessthan\id_{A\times A});b=(n\olessthan\id_A)\otimes(n\olessthan\id_A)$.
 		\begin{xalignat*}{2}
 			\lefteqn{\trans{\assoc};(\swap\otimes\id_{A\times A});(f\otimes f);f}\\
 			&= \trans{\assoc};(\swap;f\otimes f);f && \mbox{Lemma \ref{Lem:Products}(6)}\\
 			&= \trans{\assoc};(f\otimes f);f && \mbox{$f$ commutative}\\
 			&= \trans{\assoc};(f\otimes\id_{A\times A});(\id_A\otimes f);f && \mbox{Lemma \ref{Lem:Products}(6)}\\
 			&= \trans{\assoc};(f\otimes(\id_A\otimes\id_A));(\id_A\otimes f);f && \id_A\otimes\id_A=\id_{A\times A}\\
 			&= ((f\otimes\id_A)\otimes\id_A);\trans{\assoc};(\id_A\otimes f);f && \mbox{Lemma \ref{Lem:AssocSwap}(5)}\\
 			&= ((f\otimes\id_A)\otimes\id_A);(f\otimes\id_A);f && \mbox{by 1.}\\
 			&= ((f\otimes\id_A);f\otimes\id_A);f && \mbox{Lemma \ref{Lem:Products}(6)}\\
 			\lefteqn{\assoc;(\swap;\assoc\otimes\id_A);((f\otimes\id_A);f\otimes\id_A);f}\\
 			&= \assoc;(\swap;\assoc;(f\otimes\id_A);f\otimes\id_A);f && \mbox{Lemma \ref{Lem:Products}(6)}\\ 
 			&= \assoc;(\swap;(\id_A\otimes f);f\otimes\id_A);f && \mbox{$f$ associative}\\
 			&= \assoc;((f\otimes\id_A);\swap;f\otimes\id_A);f && \mbox{Lemma \ref{Lem:AssocSwap}(3)}\\
 			&= \assoc;((f\otimes\id_A);f\otimes\id_A);f && \mbox{$f$ commutative}\\ 
 			&= (f\otimes(\id_A\otimes\id_A));\assoc;(f\otimes\id_A);f && \mbox{Lemma \ref{Lem:AssocSwap}(5)}\\
 			&= (f\otimes\id_{A\times A});\assoc;(f\otimes\id_A);f && \id_A\otimes\id_A=\id_{A\times A}\\
 			&= (f\otimes\id_{A\times A});(\id_A\otimes f);f && \mbox{$f$ associative}\\
 			&= (f\otimes f);f, && \mbox{Lemma \ref{Lem:Products}(6)}
		\end{xalignat*}
		i.e., $\trans{b};(f\otimes f);f=(f\otimes f);f$. Together we conclude
		\begin{xalignat*}{2}
		 	((n\otimes n);f\olessthan f);f
		 	&= ((n\otimes n)\olessthan\id_{A\times  A});(f\otimes f);f && \mbox{Lemma \ref{Lem:Products}(6)}\\
		 	&= ((n\otimes n)\olessthan\id_{A\times  A});b;\trans{b};(f\otimes f);f && \mbox{$b$ bijective}\\
		 	&= ((n\olessthan\id_A)\otimes(n\olessthan\id_A));(f\otimes f);f && \mbox{see above}\\
		 	&= ((n\olessthan\id_A);f\otimes(n\olessthan\id_A);f);f && \mbox{Lemma \ref{Lem:Products}(6)}\\		 	
		 	&= (\Top_{A1};e\otimes\Top_{A1};e);f && \mbox{by 2.}\\
		 	&= \Top_{A1};e;(\id_A\otimes\Top_{A1};e);f && \mbox{Lemma \ref{Lem:Basics}(3)}\\
		 	&= \Top_{A1};e. && \mbox{$e$ neutral}
		\end{xalignat*}
		From 3.\ we obtain $(n\otimes n);f;n=f$.\QED
\end{enumerate}
\def\QED{}
\end{shortproof}

A relation $C:X\to X$ is called transitive iff $C;C\sqsubseteq C$, dense iff $C\sqsubseteq C;C$, asymmetric iff $C\sqcap\trans{C}=\Bot_{XX}$,
a strict-order iff $C$ is transitive and asymmetric, a linear strict-order iff $C$ is a strict-order and $\id_X\sqcup C\sqcup\trans{C}=\Top_{XX}$.

Given a strict-order $C:X\to X$ we define its associate ordering $E=\id_X\sqcup C$. It is easy to verify that $E$ is an ordering, i.e., $E$ is reflexive 
$\id_X\sqsubseteq E$, $E$ is transitive, and $E$ is antisymmetric $E\sqcap\trans{E}=\id_X$.

The next lemma verifies that a linear strict-order does always have a complement. 

\begin{lemma}\label{Lem:StrictComp}
If $C:A\to A$ is a linear strict-order and $E$ its associated ordering, then we have
\begin{enumerate}
	\item $C\sqcap\id_A=\Bot_{AA}$,
	\item $C$ is complemented with complement $\trans{E}$,
	\item $C^\star=\trans{E}$ and $E^\star=\trans{C}$,
	\item $C$ and $E$ are regular.
\end{enumerate}
\end{lemma}

\begin{shortproof}
\begin{enumerate}
	\item We have
		\begin{xalignat*}{2}
			C\sqcap\id_A
			&= (C\sqcap\id_A)\sqcap(C\sqcap\id_A)\\
			&= (C\sqcap\id_A)\sqcap\trans{(C\sqcap\id_A)} && \mbox{Lemma \ref{Lem:Basics}(1)}\\
			&= C\sqcap\id_A\sqcap\trans{C}\\
			&= \Bot_{AA}. && \mbox{$C$ linear}
		\end{xalignat*}
	\item We have $C\sqcup\trans{E}=C\sqcup\trans{C}\sqcup\id_A=\Top_{AA}$ because $C$ is linear. On the other hand, the asymmetry and (1) implies
	      $C\sqcap\trans{E}=C\sqcap(\trans{C}\sqcup\id_A)=(C\sqcap\trans{C})\sqcup(C\sqcap\id_A)=\Bot_{AA}$.
\end{enumerate}
\noindent Property 3.\ and 4.\ follow immediately from 2.\ and Lemma \ref{LemCompReg}.
\end{shortproof}

If $E$ is an order we define
\begin{xalignat*}{2}
    \ubd{E}{R}  &= \lres{\trans{R}}{E},\\
    \lbd{E}{R}  &= \lres{\trans{R}}{\trans{E}},\\
    \lub{E}{R} &= \ubd{E}{R}\sqcap\lbd{E}{\ubd{E}{R}},\\
    \glb{E}{R} &= \lbd{E}{R}\sqcap\ubd{E}{\lbd{E}{R}}.
\end{xalignat*}
For concrete relations the construction $\ubd{E}{R}$ relates $b$ with the upper bounds of the set of elements
related to $b$ in $R$, i.e., the image $\{a\mid (b,a)\in R\}$ of $b$ in $R$. Similarly, $\lbd{E}{R}$ computes the lower bounds. Finally, $\lub{E}{R}$ (resp.\ $\glb{E}{R}$)
maps $b$ to the least upper (greatest lower) bound of the image of $b$ in $R$.

\begin{lemma}\label{Lem:BoundLemma}
Let $E:A\to A$ be an ordering, and $X:B\to A$ a relation. Then we have $\ubd{E}{X};E=\ubd{E}{X}$ and $\lbd{E}{X};\trans{E}=\lbd{E}{X}$.
\end{lemma}

\begin{proof}
We only show the first property since the second can be shown analogously. We have $\trans{X};\ubd{E}{X};E=\trans{X};(\lres{\trans{X}}{E});E\sqsubseteq E;E\sqsubseteq E$ since $E$
is transitive. This implies $\ubd{E}{X};E\sqsubseteq\lres{\trans{X}}{E}=\ubd{E}{X}$. The converse inclusion follows from the reflexivity of $E$.
\end{proof}

If we consider a concrete linear strict-order and an element $a$ strictly below the least upper bound of a set $M$, then there is an element $b$ in $M$ that is already strictly greater than $a$.
We want to show this property in arbitrary Heyting categories. It turns out that we can only show this property if we use double-negation in the conclusion.

\begin{lemma}\label{Lem:DownClosed}
If $C:A\to A$ is a linear strict-order and $X:B\to A$ a relation, then $\lub{E}{X};\trans{C}\sqsubseteq(X;\trans{C})^{\star\star}$.
\end{lemma}

\begin{proof}
First of all, we have
	\begin{xalignat*}{2}
		\trans{\lub{E}{X}};(X;\trans{C})^\star
		&= \trans{\lub{E}{X}};(X;{\trans{C}}^{\star\star})^\star && \mbox{Lemma \ref{Lem:StrictComp}(4)}\\
		&= \trans{\lub{E}{X}};(\lres{\trans{X}}{{\trans{C}}^\star}) && \mbox{Lemma \ref{Lem:ResNeg}(3)}\\
		&= \trans{\lub{E}{X}};(\lres{\trans{X}}{E}) && \mbox{Lemma \ref{Lem:StrictComp}(3)}\\
		&= \trans{\lub{E}{X}};\ubd{E}{X}\\
		&\sqsubseteq \trans{(\lres{\trans{\ubd{E}{X}}}{\trans{E}})};\ubd{E}{X}\\
		&= (\rres{E}{\ubd{E}{X}});\ubd{E}{X}\\
		&\sqsubseteq E && \mbox{Lemma \ref{Lem:ResBasics}(1)}\\
		&= {\trans{C}}^\star. && \mbox{Lemma \ref{Lem:StrictComp}(3)}
	\end{xalignat*}
	This implies $\lub{E}{X};\trans{C}\sqsubseteq(X;\trans{C})^{\star\star}$ by using the Schr\"oder equivalences (Lemma \ref{Lem:Schroeder}).
\end{proof}

The following example demonstrates that we cannot remove the double-negation on the right-hand side of the inclusion of the previous lemma. It is well-known that for every Heyting algebra $L$
the collection of $L$-fuzzy relations between sets, i.e., the collection of functions $A\times B\to L$, forms a Heyting category. If both sets $A$ and $B$ are finite, we can fix a linear order on each set and represent an $L$-fuzzy relation $R$ between $A$ and $B$ by a matrix $M$ with coefficients from $L$, i.e., if $R(a,b)=x$ and $a$ is the $i$-th element of $A$ and $b$ the $j$-th element of $B$, then $M$ has an entry $x$ in row $i$ and column $j$. Composition of two relations in matrix form is based on the formula $$(Q;R)(a,c)=\bigsqcup\limits_{b\in B}Q(a,b)\sqcap R(b,c)$$ and corresponds to matrix multiplication known from linear algebra of the matrix for $Q$ and the matrix of $R$ using $\sqcap$ and $\sqcup$ instead of the multiplication and addition, respectively. Similarly, the left residual can be computed by the formula $$(\lres{Q}{R})(b,c)=\bigsqcap\limits_{a\in A}Q(a,b)\to R(a,c)$$ and corresponds to matrix multiplication of the transposed matrix for $Q$ and the matrix of $R$ using $\sqcap$ and $\to$ instead of multiplication and addition, respectively. Now consider the 3-element chain $0,u,1$ with $0$ and $1$ as smallest resp.\ greatest element and a set $A$ with one element and a set $B$ with two elements. Now we define two relations $C:B\to B$ and $X:A\to B$ in matrix form by:
$$ C=\left(\begin{array}{cc}0 & 1\\ 0 & 0\end{array}\right),~~~~X=\left(\begin{array}{cc}0 & u\end{array}\right). $$
$C$ is a linear strict-order and we have
$$ \begin{array}{rcl@{~~~~}rcl@{~~~~}rcl}
	\ubd{E}{X}&=&\left(\begin{array}{cc}0 & 1\end{array}\right), & \lbd{E}{\ubd{E}{X}}&=&\left(\begin{array}{cc}1 & 1\end{array}\right), & \lub{E}{X}&=&\left(\begin{array}{cc}0 & 1\end{array}\right),\\
	\lub{E}{X};\trans{C}&=&\left(\begin{array}{cc}1 & 0\end{array}\right), & X;\trans{C}&=&\left(\begin{array}{cc}u & 0\end{array}\right), & (X;\trans{C})^{\star\star}&=&\left(\begin{array}{cc}1 & 0\end{array}\right).
\end{array} $$

The relation algebraic version of a power set is given by a so-called relational (or direct) power.

\begin{definition}
An object ${\mathcal P}(A)$ together with a relation $\epsi:A\to{\mathcal P}(A)$ is called
a relational (or direct) power of $A$ iff
$$ \syq{\epsi}{\epsi}=\id_{{\mathcal P}(A)}~~~~\mbox{and}~~~~\syq{Q}{\epsi}\mbox{ is total for every }Q:A\to B.$$
\end{definition}

Please note that $\syq{\trans{R}}{\epsi}$ is a map for every relation $R:B\to A$. In fact, this construction is the existential image of $R$, i.e., $x$ is mapped by $\syq{\trans{R}}{\epsi}$ to the set
$\{y\mid (x,y)\in R\}$ for concrete relations. Furthermore, we have the following \cite{Schmidt,RedBible,Winter,Book}.

\begin{lemma}\label{Lem:PowerBasic}
Suppose $R:A\to B$ is a relation. Then we have $\syq{\trans{R}}{\epsi};\trans{\epsi}=R$.
\end{lemma}

The fourth axiom of Tarski's axiomatization of the real numbers requires for all non-empty subsets $X,Y$ of the real numbers with $x<y$ for every $x\in X$ and $y\in Y$ the existence of an element $z$ with $x\leq z$ and $z\leq y$ for all $x\in X$ and $y\in Y$. Our original definition \cite{PartI} of a real number object used the relational power from above which includes the empty set in the case of concrete relations. This is not correct but, fortunately, the results presented in \cite{PartI} did not rely on Axiom 4, i.e., Axiom 4 was never used. In the current paper we fix this mistake and switch to the non-empty relational power that corresponds to the set of non-empty subsets in the case of concrete relations.

\begin{definition}
An object $\nepower{A}$ together with a relation $\epsilon:A\to\nepower{A}$ is called
a non-empty relational power of $A$ iff
$$ \syq{\epsilon}{\epsilon}=\id_{\nepower{A}}~~~~\mbox{and}~~~~\dom{\syq{Q}{\epsilon}}=\dom{\trans{Q}}\mbox{ for every }Q:A\to B.$$
\end{definition}

We want to show that non-empty relational powers exists if relational powers and splittings exists. Given a partial equivalence relation $X:A\to A$, i.e., $X$ is symmetric ($\trans{X}=X$) and transitive, 
an object $B$ together with a relation $R:B\to A$ is called a splitting of $X$ iff $R;\trans{R}=\id_B$ and $\trans{R};R=X$. Intuitively, the object $B$ consists of all (existing) equivalence classes of $X$ 
and $R$ relates such an equivalence class with its elements. Please note that requiring the existence of splittings is not really an additional assumption since every Heyting category can be fully embedded 
into a Heyting category with all splittings \cite{Freyd,Winter}.

\begin{theorem}
If $i:C\to{\mathcal P}(A)$ splits $\dom{\trans{\epsi}}$, then $C$ together with $\epsi;\trans{i}$ is a non-empty relational power.
\end{theorem}

\begin{shortproof}
First of all, we have
\begin{xalignat*}{2}
	\syq{\epsi;\trans{i}}{\epsi;\trans{i}}
	&= i;\syq{\epsi}{\epsi};\trans{i} && \mbox{Lemma \ref{Lem:ResBasics}(4)}\\
	&= i;\trans{i} && \mbox{definition $\epsi$}\\
	&= \id_C. && \mbox{definition $i$}
\end{xalignat*}
Now, suppose $Q:A\to B$. Then we have
\begin{xalignat*}{2}
	\lefteqn{\dom{\trans{Q}}}\\
	&= \id_B\sqcap\trans{Q};Q\\
	&= \id_B\sqcap\syq{Q}{\epsi};\trans{\epsi};\epsi;\syq{\epsi}{Q} && \mbox{Lemma \ref{Lem:PowerBasic}}\\
	&= \id_B\sqcap\syq{Q}{\epsi};\trans{\epsi};\epsi;\trans{\syq{Q}{\epsi}} && \mbox{Lemma \ref{Lem:ResBasics}(3)}\\
	&= \id_B\sqcap\syq{Q}{\epsi};\trans{\syq{Q}{\epsi}}\sqcap\syq{Q}{\epsi};\trans{\epsi};\epsi;\trans{\syq{Q}{\epsi}} && \mbox{$\syq{Q}{\epsi}$ total}\\
	&= \id_B\sqcap\syq{Q}{\epsi};(\id\sqcap\trans{\epsi};\epsi);\trans{\syq{Q}{\epsi}} && \mbox{Lemma \ref{Lem:Maps}(3)}\\
	&= \id_B\sqcap\syq{Q}{\epsi};\trans{i};i;\trans{\syq{Q}{\epsi}} && \mbox{definition $i$}\\
	&= \id_B\sqcap\syq{Q}{\epsi;\trans{i}};\trans{\syq{Q}{\epsi;\trans{i}}} && \mbox{Lemma \ref{Lem:ResBasics}(3)}\\
	&= \dom{\syq{Q}{\epsi;\trans{i}}}. 
\end{xalignat*}
\vspace{-5ex}

{~}
\end{shortproof}

Last but not least, we want to show Lemma \ref{Lem:PowerBasic} also for non-empty powers.

\begin{lemma}\label{Lem:TPowerBasic}
Suppose $R:A\to B$ is a relation. Then we have $\syq{\trans{R}}{\epsilon};\trans{\epsilon}=R$.
\end{lemma}

\begin{proof}
First of all, we have $\syq{\trans{R}}{\epsilon};\trans{\epsilon}\sqsubseteq(\rres{R}{\trans{\epsilon}}) ;\trans{\epsilon}\sqsubseteq R$.
\begin{xalignat*}{2}
	R
	&= \dom{R};R\\
	&= \dom{\syq{\trans{R}}{\epsilon}};R && \mbox{definition $\epsilon$}\\
	&\sqsubseteq \syq{\trans{R}}{\epsilon};\trans{\syq{\trans{R}}{\epsilon}};R\\
	&= \syq{\trans{R}}{\epsilon};\syq{\epsilon}{\trans{R}};R && \mbox{Lemma \ref{Lem:ResBasics}(3)}\\	
	&\sqsubseteq \syq{\trans{R}}{\epsilon};(\rres{\trans{\epsilon}}{R});R\\
	&\sqsubseteq \syq{\trans{R}}{\epsilon};\trans{\epsilon}.\tag*{~}
\end{xalignat*}
\vspace{-4ex}

{~}
\end{proof}

\section{Real Number Object}

\vspace{-1ex}

In this section we want to recall the results of \cite{PartI}. We start with Tarski's axioms as they were stated in \cite{Tarski}.
His axioms are based on the language ${\mathbb R},<,+,1$:

\begin{description}
    \item[{\rm Axiom 1}:] If $x\ne y$, then $x<y$ or $y<x$.
    \item[{\rm Axiom 2}:] If $x<y$, then $y\nless x$.
    \item[{\rm Axiom 3}:] If $x<z$, then there is a $y$ such that $x<y$ and $y<z$.
    \item[{\rm Axiom 4}:] If $\emptyset\ne X\subseteq{\mathbb R}$ and $\emptyset\ne Y\subseteq{\mathbb R}$ so that for
    every $x\in X$ and every $y\in Y$ we have $x<y$, then there is a $z$ so that for all $x\in X$
    and $y\in Y$ we have $x\leq z$ and $z\leq y$ ($x\leq y$ shorthand for $x<y$ or $x=y$).
    \item[{\rm Axiom 5}:] $x+(y+z)=(x+z)+y$.
    \item[{\rm Axiom 6}:] For every $x$ and $y$ there is a $z$ such that $x=y+z$.
    \item[{\rm Axiom 7}:] If $x+z<y+t$, then $x<y$ or $z<t$.
    \item[{\rm Axiom 8}:] $1\in{\mathbb R}$.
    \item[{\rm Axiom 9}:] $1<1+1$.
\end{description}

\vspace{-1ex}

The axioms above can be translated into the language of relations leading to the definition below. Please note that we added Axiom 0 that states that $\add$ is a map explicitly since we are dealing with relations
rather than functions.

\vspace{-1ex}

\begin{definition}\label{Def:RealNumberObject}
An object ${\mathbb R}$ together with three relations $\one:1\to{\mathbb R}$, $C:{\mathbb R}\to{\mathbb R}$ and
$\add:{\mathbb R}\times {\mathbb R}\to{\mathbb R}$ is called a real number object if the following holds:
\begin{enumerate}\setcounter{enumi}{-1}
    \item $\add$ is a map.
    \item $\id_{\mathbb R}\sqcup C\sqcup\trans{C}=\Top_{\mathbb{R R}}$.
    \item $C\sqcap\trans{C}=\Bot_{\mathbb{R R}}$.
    \item $C\sqsubseteq C;C$.
    \item $\lres{\epsilon}{(\rres{C}{\trans{\epsilon}})}\sqsubseteq(\lres{\epsilon}{(C\sqcup\id_{\mathbb R})});\trans{(\lres{\epsilon}{\trans{(C\sqcup\id_{\mathbb R})}})}$.
    \item $(\id_{\mathbb R}\otimes\add);\add=(\id_{\mathbb R}\otimes\swap);\assoc;(\add\otimes\id_{\mathbb R});\add$.
    \item $\trans{\pi};\add=\Top_{\mathbb{R R}}$.
    \item $\add;C;\trans{\add}\sqsubseteq\pi;C;\trans{\pi}\sqcup\rho;C;\trans{\rho}$.
    \item $\one$ is a map, i.e., a point.
    \item $\one\sqsubseteq \one;(\id_{\mathbb R}\olessthan\id_{\mathbb R});\add;\trans{C}$.
\end{enumerate}
\end{definition}

First we define abstract versions of the number $0$ and of the negation operation on the real numbers by
$\zero=\Top_{1{\mathbb R}};(\trans{\add}\sqcap\trans{\pi});\rho$ and $\lbdnus=\trans{\pi};(\add;\trans{Z}\sqcap\rho)$.

The first main result of \cite{PartI} is the following theorem.

\begin{theorem}\label{Th:AddGroup}
The quadruple $({\mathbb R},\zero,\add,\lbdnus)$ is an abelian group.
\end{theorem}

The second main result is concerned with the strict-order $C$.

\begin{theorem}
The relation $C:{\mathbb R}\to {\mathbb R}$ is a dense strict linear order.
\end{theorem}

Last but not least, the final result of \cite{PartI} addresses the monotonicity of $\add$.

\begin{theorem}\label{Th:AddMono}
We have the following
\begin{enumerate}
	\item $\add$ is strictly monotone in each parameter, i.e., $(\id_{\mathbb R}\otimes C);\add\sqsubseteq\add;C$ and $(C\otimes\id_{\mathbb R});\add\sqsubseteq\add;C$,
	\item $\add$ is strictly monotone. i.e., $(C\otimes C);\add\sqsubseteq\add;C$,
	\item $\add$ is monotone, i.e., $(E\otimes E);\add\sqsubseteq\add;E$.
\end{enumerate}
\end{theorem}

\section{Least-Upper-Bound Property}

In this final section we want to show the least-upper-bound property of a real number object. But first we will show that adding a constant to a number and subtracting the same constant are inverse operations. 

\begin{lemma}
Suppose $p:1\to{\mathbb R}$ is a point. Then $(\id_{\mathbb R}\olessthan \Top_{{\mathbb R}1};p);\add$ is strictly monotone and a bijective map with 
$\trans{((\id_{\mathbb R}\olessthan \Top_{{\mathbb R}1};p);\add)}=(\id_{\mathbb R}\olessthan \Top_{{\mathbb R}1};p;\lbdnus);\add$.
\end{lemma}

\begin{proof}
First of all, $(\id_{\mathbb R}\olessthan \Top_{{\mathbb R}1};p);\add$ is a map because $p$ and $\add$ are. Now we show $(\id_{\mathbb R}\olessthan \Top_{{\mathbb R}1};p);\add;(\id_{\mathbb R}\olessthan \Top_{{\mathbb R}1};p;\lbdnus);\add=\id_{\mathbb R}$ by
computing
\begin{xalignat*}{2}
	\lefteqn{(\id_{\mathbb R}\olessthan \Top_{{\mathbb R}1};p);\add;(\id_{\mathbb R}\olessthan \Top_{{\mathbb R}1};p;\lbdnus);\add}\\
	&= ((\id_{\mathbb R}\olessthan \Top_{{\mathbb R}1};p);\add\olessthan \Top_{{\mathbb R}1};p;\lbdnus);\add && \mbox{Lemma \ref{Lem:Basics}(3)}\\
	&= ((\id_{\mathbb R}\olessthan \Top_{{\mathbb R}1};p)\olessthan \Top_{{\mathbb R}1};p;\lbdnus);(\add\otimes\id_{\mathbb R});\add && \mbox{Lemma \ref{Lem:Products}(6)}\\
	&= ((\id_{\mathbb R}\olessthan \Top_{{\mathbb R}1};p)\olessthan \Top_{{\mathbb R}1};p;\lbdnus);\trans{\assoc};(\id_{\mathbb R}\otimes\add);\add && \mbox{Lemma \ref{Lem:GroupProps}(1)}\\
	&= ((\id_{\mathbb R}\olessthan (\Top_{{\mathbb R}1};p\olessthan \Top_{{\mathbb R}1};p;\lbdnus));(\id_{\mathbb R}\otimes\add);\add && \mbox{Lemma \ref{Lem:AssocSwap}(3)}\\
	&= ((\id_{\mathbb R}\olessthan \Top_{{\mathbb R}1};p;(\id_{\mathbb R}\olessthan\lbdnus));(\id_{\mathbb R}\otimes\add);\add && \mbox{Lemma \ref{Lem:Maps}(3)}\\
	&= ((\id_{\mathbb R}\olessthan \Top_{{\mathbb R}1};p;(\id_{\mathbb R}\olessthan\lbdnus);\add);\add && \mbox{Lemma \ref{Lem:Products}(6)}\\
	&= ((\id_{\mathbb R}\olessthan \Top_{{\mathbb R}1};p;\Top_{{\mathbb R}1};0);\add && \mbox{Theorem \ref{Th:AddGroup}}\\
	&= ((\id_{\mathbb R}\olessthan \Top_{{\mathbb R}1};\Top_{11};0);\add && \mbox{$p$ total}\\
	&= ((\id_{\mathbb R}\olessthan \Top_{{\mathbb R}1};0);\add && \mbox{$\Top_{11}=\id_1$}\\
	&= \id_{\mathbb R}. && \mbox{Theorem \ref{Th:AddGroup}}
\end{xalignat*}
The property $(\id_{\mathbb R}\olessthan \Top_{{\mathbb R}1};p;\lbdnus);\add;(\id_{\mathbb R}\olessthan \Top_{{\mathbb R}1};p);\add=\id_{\mathbb R}$ can be shown analogously.
Using Lemma \ref{Lem:Basics}(1) the two properties above imply $(\id_{\mathbb R}\olessthan \Top_{{\mathbb R}1};p);\add\sqsubseteq\trans{((\id_{\mathbb R}\olessthan \Top_{{\mathbb R}1};p;\lbdnus);\add)}$ and
$(\id_{\mathbb R}\olessthan \Top_{{\mathbb R}1};p;\lbdnus);\add\sqsubseteq\trans{((\id_{\mathbb R}\olessthan \Top_{{\mathbb R}1};p);\add)}$, and, hence,
$\trans{((\id_{\mathbb R}\olessthan \Top_{{\mathbb R}1};p);\add)}=(\id_{\mathbb R}\olessthan \Top_{{\mathbb R}1};p;\lbdnus);\add$ and that $((\id_{\mathbb R}\olessthan \Top_{{\mathbb R}1};p);\add$ is bijective.
It remains to show that $(\id_{\mathbb R}\olessthan \Top_{{\mathbb R}1};p);\add$ is strictly monotone. The computation
\begin{xalignat*}{2}
  	C;(\id_{\mathbb R}\olessthan \Top_{{\mathbb R}1};p);\add
  	&= (C\olessthan \Top_{{\mathbb R}1};p);\add && \mbox{Lemma \ref{Lem:Basics}(3)}\\
  	&= (\id_{\mathbb R}\olessthan \Top_{{\mathbb R}1};p);(C\otimes\id_{\mathbb R});\add && \mbox{Lemma \ref{Lem:Products}(6)}\\
  	&\sqsubseteq (\id_{\mathbb R}\olessthan \Top_{{\mathbb R}1};p);\add;C. && \mbox{Theorem \ref{Th:AddMono}(1)}
\end{xalignat*}
shows that property.
\end{proof}

We now use $\one$ for $p$ and define $\su=(\id_{\mathbb R}\olessthan \Top_{{\mathbb R}1};\one);\add$ and $\pr=(\id_{\mathbb R}\olessthan\Top_{{\mathbb R}1};\one;\lbdnus);\add$. From
the previous lemma we obtain that $\su$ as well as $\pr$ are monotone, strictly monotone, bijective, and $\trans{\su}=\pr$.

\begin{lemma}\label{Lem:Props}
We have
\begin{enumerate}
      \item $\su;\lbdnus=\lbdnus;\pr$,
	\item $0\sqsubseteq\one;\trans{C}$,
	\item $\su\sqsubseteq C$ and $\pr\sqsubseteq\trans{C}$,
	\item $C$ is total and surjective.
\end{enumerate}
\end{lemma}

\begin{shortproof}
\begin{enumerate}
	\item From the computation
		\begin{xalignat*}{2}
			\lbdnus;\pr;\lbdnus
			&= \lbdnus;(\id_{\mathbb R}\olessthan\Top_{{\mathbb R}1};\one;\lbdnus);\add;\lbdnus\\
			&= (\lbdnus\olessthan\Top_{{\mathbb R}1};\one;\lbdnus);\add;\lbdnus && \mbox{Lemma \ref{Lem:Basics}(3)}\\
			&= (\id_{\mathbb R}\olessthan\Top_{{\mathbb R}1};\one);(\lbdnus\otimes\lbdnus);\add;\lbdnus && \mbox{Lemma \ref{Lem:Products}(6)}\\
			&= (\id_{\mathbb R}\olessthan\Top_{{\mathbb R}1};\one);\add && \mbox{Lemma \ref{Lem:GroupProps}(4)}\\
			&= \su
	    	\end{xalignat*}
	    	we immediately obtain the assertion since $\lbdnus$ is a bijective map.
	\item First of all, from the fact that $\one$ is total we obtain $\one;\Top_{{\mathbb R}1}=\Top_{11}=\id_1$. Now we compute
		\begin{xalignat*}{2}
			0
			&= \one;\Top_{{\mathbb R}1};0 && \mbox{see above}\\
			&= \one;(\id_{\mathbb R}\olessthan\lbdnus);\add && \mbox{Theorem \ref{Th:AddGroup}}\\
			&= (\one\olessthan\one;\lbdnus);\add && \mbox{Lemma \ref{Lem:Maps}(3)}\\
			&= (\one\olessthan\Top_{11};\one;\lbdnus);\add && \Top_{11}=\id_1\\	
			&= \one;(\id_{\mathbb R}\olessthan\Top_{{\mathbb R}1};\one;\lbdnus);\add &&  \mbox{Lemma \ref{Lem:Basics}(3)}\\
			&= \one;\pr\\
			&\sqsubseteq \one;(\id_{\mathbb R}\olessthan\id_{\mathbb R});\add;\trans{C};\pr && \mbox{Axiom 9}\\
			&= (\one\olessthan\one);\add;\trans{C};\pr && \mbox{Lemma \ref{Lem:Maps}(3)}\\
			&= (\one\olessthan\Top_{11}\one);\add;\trans{C};\pr && \Top_{11}=\id_1\\
			&= \one;(\id_{\mathbb R}\olessthan\Top_{{\mathbb R}1}\one);\add;\trans{C};\pr && \mbox{Lemma \ref{Lem:Basics}(3)}\\
			&= \one;\su;\trans{C};\pr\\
			&= \one;\su;\pr;\trans{C} && \mbox{$\pr$ strictly monotone}\\
			&= \one;\trans{C}. && \mbox{$\su$ and $\pr$ inverse}
		\end{xalignat*}
	\item We obtain
		\begin{xalignat*}{2}
			\id_{\mathbb R}
			&= (\id_{\mathbb R}\olessthan\Top_{{\mathbb R}1};0);\add && \mbox{Theorem \ref{Th:AddGroup}}\\
			&\sqsubseteq (\id_{\mathbb R}\olessthan\Top_{{\mathbb R}1};\one;\trans{C});\add && \mbox{by 2.}\\
			&= (\id_{\mathbb R}\olessthan\Top_{{\mathbb R}1};\one);(\id_{\mathbb R}\otimes\trans{C});\add && \mbox{Lemma \ref{Lem:Products}(6)}\\	
			&\sqsubseteq (\id_{\mathbb R}\olessthan\Top_{{\mathbb R}1};\one);\add;\trans{C} && \mbox{$\add$ strictly monotone}\\
			&=\su;\trans{C}
		\end{xalignat*}
		from which we conclude $\su\sqsubseteq C;\trans{\su};\su\sqsubseteq C$ since $\su$ is univalent. The second inclusion follows from the first by
		$\pr=\trans{\su}\sqsubseteq\trans{C}$.
	\item Both properties follow immediately form 3.\ because $\su$and $\pr$ are total.
\end{enumerate}
\end{shortproof}

The following lemma will be needed in the proof of the least-upper-bound property.

\begin{lemma}\label{Lem:XCProps}
Suppose $X:A\to{\mathbb R}$. Then we have
\begin{enumerate}
	\item $\dom{X}=\dom{X;\trans{C}}$,	
	\item $\ubd{E}{X}=\ubd{E}{X;\trans{C}}$.
\end{enumerate}
\end{lemma}

\begin{shortproof}
\begin{enumerate}
	\item The inclusion $\sqsupseteq$ follows from
		\begin{xalignat*}{2}
			\dom{X;\trans{C}}
			&= \id_{\mathbb R}\sqcap X;\trans{C};C;\trans{X}\\
			&= \id_{\mathbb R}\sqcap \dom{X};X;\trans{C};C;\trans{X}\\
			&\sqsubseteq \dom{X};(\trans{\dom{X}}\sqcap X;\trans{C};C;\trans{X})\\
			&\sqsubseteq \dom{X};\trans{\dom{X}}\\
			&= \dom{X}, && \mbox{Lemma \ref{Lem:Basics}(1\&2)}
		\end{xalignat*}
		and the opposite inclusion from
		\begin{xalignat*}{2}
			\dom{X}
			&= \id_{\mathbb R}\sqcap X;\trans{X}\\
			&= \id_{\mathbb R}\sqcap (X\sqcap X;\trans{C};C);\trans{X} && \mbox{Lemma \ref{Lem:Props}(4)}\\
			&\sqsubseteq \id_{\mathbb R}\sqcap (\id_{\mathbb R}\sqcap X;\trans{C};C;\trans{X});X;\trans{X}\\	
			&= \id_{\mathbb R}\sqcap\dom{X;\trans{C}};X;\trans{X}\\
			&\sqsubseteq \dom{X;\trans{C}};(\trans{\dom{X;\trans{C}}}\sqcap X;\trans{X})\\
			&\sqsubseteq \dom{X;\trans{C}};\trans{\dom{X;\trans{C}}}\\			
			&= \dom{X;\trans{C}}. && \mbox{Lemma \ref{Lem:Basics}(1\&2)}
		\end{xalignat*}
	\item First of all, we have
		\begin{xalignat*}{2}
			\trans{X};\ubd{E}{X;\trans{C}}\sqcap\trans{C}
			&\sqsubseteq \trans{X};\ubd{E}{X;\trans{C}}\sqcap\trans{C};\trans{C} && \mbox{$C$ transitive}\\\
			&\sqsubseteq \trans{C};(C;\trans{X};\ubd{E}{X;\trans{C}}\sqcap\trans{C})\\
			&= \trans{C};(C;\trans{X};(\lres{C;\trans{X}}{E})\sqcap\trans{C})\\
			&\sqsubseteq \trans{C};(E\sqcap\trans{C}) && \mbox{Lemma \ref{Lem:ResBasics}(1)}\\
			&= \trans{C};\Bot_{{\mathbb R}{\mathbb R}} && \mbox{Lemma \ref{Lem:StrictComp}(2)}\\
			&= \Bot_{{\mathbb R}{\mathbb R}}.
		\end{xalignat*}
		This implies
		\begin{xalignat*}{2}
			\lefteqn{\trans{X};\ubd{E}{X;\trans{C}}}\\
			&= \trans{X};\ubd{E}{X;\trans{C}}\sqcap(\trans{C}\sqcup E) && \mbox{$C$ linear}\\
			&= (\trans{X};\ubd{E}{X;\trans{C}}\sqcap\trans{C})\sqcup(\trans{X};\ubd{E}{X;\trans{C}}\sqcap E)\\
			&= \trans{X};\ubd{E}{X;\trans{C}}\sqcap E, && \mbox{see above}
            	\end{xalignat*}
            	i.e., $\trans{X};\ubd{E}{X;\trans{C}}\sqsubseteq E$. We conclude $\ubd{E}{X;\trans{C}}\sqsubseteq\lres{\trans{X}}{E}=\ubd{E}{X}$. For the opposite inclusion consider
            	\begin{xalignat*}{2}
            		C;\trans{X};\ubd{E}{X}
            		&= C;\trans{X};(\lres{\trans{X}}{E})\\	
            		&\sqsubseteq C;E && \mbox{Lemma \ref{Lem:ResBasics}(1)}\\
            		&= C;C\sqcup C\\
            		&= C && \mbox{$C$ transitive}\\
            		&\sqsubseteq E
            	\end{xalignat*}
            	from which we conclude $\ubd{E}{X}\sqsubseteq\lres{C;\trans{X}}{E}=\ubd{E}{X;\trans{C}}$.
\end{enumerate}
\end{shortproof}

Now, we are ready to show the least-upper-bound property. A relation $X:A\to {\mathbb R}$ can be seen a a collection of subsets of ${\mathbb R}$ indexed by $A$, i.e., every $a\in A$ is related to its image under $X$.
The element $a$ is in the domain of $X$ iff its image is non-empty. Therefore, the relation $\dom{X}\sqcap\dom{\ubd{E}{X}}$ relates an element $a$ to itself iff its image and the upper bound of its image are not empty.
The least-upper-bound property now states that least upper bound for such a set exists, i.e., that $\dom{X}\sqcap\dom{\ubd{E}{X}}\sqsubseteq\lub{E}{X}$. This is our main theorem of the paper.	

\begin{theorem}[Least-Upper-Bound Property]\label{Th:LUP}
For every relation $X:A\to{\mathbb R}$ we have $\dom{X}\sqcap\dom{\ubd{E}{X}}\sqsubseteq\dom{\lub{E}{X}}$.
\end{theorem}

\begin{proof}
First of all, we have
\begin{xalignat*}{2}
	\lefteqn{\epsi;\trans{\syq{C;\trans{X}}{\epsilon}};\syq{\trans{\ubd{E}{X}}}{\epsilon});\trans{\epsilon}}\\
	&= \trans{(\syq{C;\trans{X}}{\epsilon};\trans{\epsilon})}\syq{\trans{\ubd{E}{X}}}{\epsilon});\trans{\epsilon}\\
	&= C;\trans{X};\ubd{E}{X} && \mbox{Lemma \ref{Lem:TPowerBasic}}\\
	&= C;\trans{X};(\lres{\trans{X}}{E})\\
	&\sqsubseteq C;E && \mbox{Lemma \ref{Lem:ResBasics}(1)}\\
	&= C;(C\sqcup\id_{\mathbb R})\\
	&= C;C\sqcup C\\
	&= C && \mbox{$C$ transitive}
\end{xalignat*}
which immediately implies $\trans{\syq{C;\trans{X}}{\epsilon}};\syq{\trans{\ubd{E}{X}}}{\epsilon}\sqsubseteq\lres{\epsi}{(\rres{C}{\trans{\epsi}})}$. We obtain
\begin{xalignat*}{2}
	\lefteqn{\dom{X}\sqcap\dom{\ubd{E}{X}}}\\
	&= \dom{X};\dom{\ubd{E}{X}} && \mbox{Lemma \ref{Lem:Basics}(2)}\\	
	&= \dom{X;\trans{C}};\dom{\ubd{E}{X}} && \mbox{Lemma \ref{Lem:XCProps}(1)}\\	
	&= \dom{\syq{C;\trans{X}}{\epsilon}};\dom{\syq{\trans{\ubd{E}{X}}}{\epsilon}} && \mbox{Definition $\epsilon$}\\	
	&\sqsubseteq \syq{C;\trans{X}}{\epsilon};\trans{\syq{C;\trans{X}}{\epsilon}};\syq{\trans{\ubd{E}{X}}}{\epsilon};\trans{\syq{\trans{\ubd{E}{X}}}{\epsilon}}\span\omit\span\omit\\	
	&\sqsubseteq \syq{C;\trans{X}}{\epsilon};(\lres{\epsilon}{(\rres{C}{\trans{\epsilon}})});\trans{\syq{\trans{\ubd{E}{X}}}{\epsilon}} && \mbox{see above}\\
	&\sqsubseteq \syq{C;\trans{X}}{\epsilon};(\lres{\epsilon}{E});\trans{(\lres{\epsilon}{\trans{E})}};\trans{\syq{\trans{\ubd{E}{X}}}{\epsilon}} && \mbox{Axiom 4}\\
	&\sqsubseteq (\lres{\epsilon;\trans{\syq{C;\trans{X}}{\epsilon}}}{E});\trans{(\lres{\epsilon;\trans{\syq{\trans{\ubd{E}{X}}}{\epsilon}}}{\trans{E})}} && \mbox{Lemma \ref{Lem:ResBasics}(2)}\\
	&= (\lres{\trans{(X;\trans{C})}}{E});\trans{(\lres{\trans{\ubd{E}{X}}}{\trans{E})}} && \mbox{Lemma \ref{Lem:TPowerBasic}}\\
	&= \ubd{E}{X;\trans{C}};\trans{\lbd{E}{\ubd{E}{X}}} && \mbox{Definition}\\
	&= \ubd{E}{X};\trans{\lbd{E}{\ubd{E}{X}}}. && \mbox{Lemma \ref{Lem:XCProps}(2)}
\end{xalignat*}
This immediately implies
\begin{xalignat*}{2}
	\lefteqn{\dom{X}\sqcap\dom{\ubd{E}{X}}}\\
	&= \dom{X}\sqcap\dom{\ubd{E}{X}}\sqcap\trans{(\dom{X}\sqcap\dom{\ubd{E}{X}})} && \mbox{Lemma \ref{Lem:Basics}(1)}\\
	&\sqsubseteq \ubd{E}{X};\trans{\lbd{E}{\ubd{E}{X}}}\sqcap\lbd{E}{\ubd{E}{X}};\trans{\ubd{E}{X}} && \mbox{see above}\\
	&\sqsubseteq (\ubd{E}{X}\sqcap\lbd{E}{\ubd{E}{X}};\trans{\ubd{E}{X}};\lbd{E}{\ubd{E}{X}});\trans{\lbd{E}{\ubd{E}{X}}}\span\omit\span\omit\\
	&= (\ubd{E}{X}\sqcap\lbd{E}{\ubd{E}{X}};\trans{\ubd{E}{X}};(\lres{\trans{\ubd{E}{X}}}{\trans{E}}));\trans{\lbd{E}{\ubd{E}{X}}}\span\omit\span\omit\\
	&\sqsubseteq (\ubd{E}{X}\sqcap\lbd{E}{\ubd{E}{X}};\trans{E});\trans{\lbd{E}{\ubd{E}{X}}} && \mbox{Lemma \ref{Lem:ResBasics}(1)}\\
	&= (\ubd{E}{X}\sqcap\lbd{E}{\ubd{E}{X}});\trans{\lbd{E}{\ubd{E}{X}}} && \mbox{Lemma \ref{Lem:BoundLemma}}\\
	&= \lub{E}{X};\trans{\lbd{E}{\ubd{E}{X}}}
\end{xalignat*}
and
\begin{xalignat*}{2}
	\lefteqn{\dom{X}\sqcap\dom{\ubd{E}{X}}}\\
	&= \dom{X}\sqcap\dom{\ubd{E}{X}}\sqcap\trans{(\dom{X}\sqcap\dom{\ubd{E}{X}})} && \mbox{Lemma \ref{Lem:Basics}(1)}\\
	&\sqsubseteq \ubd{E}{X};\trans{\lbd{E}{\ubd{E}{X}}}\sqcap\lbd{E}{\ubd{E}{X}};\trans{\ubd{E}{X}} && \mbox{see above}\\
	&\sqsubseteq (\ubd{E}{X};\trans{\lbd{E}{\ubd{E}{X}}};\ubd{E}{X}\sqcap\lbd{E}{\ubd{E}{X}});\trans{\ubd{E}{X}}\span\omit\span\omit\\
	&= (\ubd{E}{X};\trans{(\trans{\ubd{E}{X}};(\lres{\trans{\ubd{E}{X}}}{\trans{E}}))}\sqcap\lbd{E}{\ubd{E}{X}});\trans{\ubd{E}{X}}\span\omit\span\omit\\
	&\sqsubseteq (\ubd{E}{X};E\sqcap\lbd{E}{\ubd{E}{X}});\trans{\ubd{E}{X}} && \mbox{Lemma \ref{Lem:ResBasics}(1)}\\
	&= (\ubd{E}{X}\sqcap\lbd{E}{\ubd{E}{X}});\trans{\ubd{E}{X}} && \mbox{Lemma \ref{Lem:BoundLemma}}\\
	&= \lub{E}{X};\trans{\ubd{E}{X}}.
\end{xalignat*}
Together we obtain
\begin{xalignat*}{2}
	\lefteqn{\dom{X}\sqcap\dom{\ubd{E}{X}}}\\
	&\sqsubseteq\lub{E}{X};\trans{\ubd{E}{X}}\sqcap\lub{E}{X};\trans{\lbd{E}{\ubd{E}{X}}} && \mbox{see above}\\
	&= \lub{E}{X};(\trans{\ubd{E}{X}}\sqcap\trans{\lbd{E}{\ubd{E}{X}}}) && \mbox{Lemma \ref{Lem:Maps}(3)}\\
	&= \lub{E}{X};\trans{\lub{E}{X}}
\end{xalignat*}
which immediately implies $$\dom{X}\sqcap\dom{\ubd{E}{X}}\sqsubseteq\dom{\lub{E}{X}}$$ since $\dom{X}\sqcap\dom{\ubd{E}{X}}\sqsubseteq\id_A$.
\end{proof}

One would expect that the inclusion of the previous theorem is in fact an equation. However, we are only able to show this for arbitrary Heyting categories if we require an additional regularity condition.

\begin{lemma}
If $X:A\to{\mathbb R}$ with $X;\trans{C}$ regular, then we have $\dom{X}\sqcap\dom{\ubd{E}{X}}=\dom{\lub{E}{X}}$.
\end{lemma}

\begin{proof}
By Theorem \ref{Th:LUP} it is sufficient to show $\dom{\lub{E}{X}}\sqsubseteq\dom{X}\sqcap\dom{\ubd{E}{X}}$. First of all, we have obviously have
$\dom{\lub{E}{X}}\sqsubseteq\dom{\ubd{E}{X}}$. Furthermore we have
\begin{xalignat*}{2}
	\dom{\lub{E}{X}}
	&= \dom{\lub{E}{X};\trans{C}} && \mbox{Lemma \ref{Lem:XCProps}(1)}\\
	&\sqsubseteq \dom{(X;\trans{C})^{\star\star}} && \mbox{Lemma \ref{Lem:DownClosed}}\\
	&= \dom{X;\trans{C}} && \mbox{assumption}\\
	&= \dom{X}. && \mbox{Lemma \ref{Lem:XCProps}(1)}
\end{xalignat*}
\end{proof}

\section{Conclusion and Future Work}

The current paper has shown the least-upper-bound property for a real number object in a Heyting category. This is the first step for showing that this additive group is Archimedean. 
For a next step one first has to define the the operation of summing up $n$ copies of an element $a$, i.e., a map ${\mathbb N}\times{\mathbb R}\to{\mathbb R}$. This requires 
either an external object of the natural numbers or to identify the natural numbers within the real number object.

Another paper will concentrate on the multiplicative group of a real number object. The definition of the multiplication operation requires the Archimedean property
and shows that the multiplication of natural number has a unique extension in the real numbers.

Last but not least, we would like to study the topology induced by the order structure on a real number object using the relation algebraic approach to topological spaces \cite{SchmidtWinter}.


\begin{thebibliography}{99}

\bibitem{Freyd}
Freyd P., Scedrov A.: Categories, allegories.
North-Holland Mathematical Library Vol. 39,
North-Holland, Amsterdam (1990).

\bibitem{oli80}
Olivier J.P., Serrato D.: Cat{\'e}gories de Dedekind. Morphismes
dans les Cat{\'e}gories de Schr\"oder. C.R.\ Acad.\ Sci.\ Paris
290 (1980) 939-941.

\bibitem{oli95}
Olivier J.P., Serrato D.: Squares and Rectangles in Relational
Categories - Three Cases: Semilattice, Distributive lattice and
Boolean Non-unitary. Fuzzy sets and systems 72 (1995), 167-178.

\bibitem{Schmidt}
  Schmidt G., Str\"ohlein T.:
  Relations and graphs.
   Discrete mathematics for computer scientists.
 EATCS Monographs on Theoretical Computer Science,
  Springer, Berlin (1993).

\bibitem{RedBible}
  Schmidt G.:
  Relational mathematics.
 Encyplopedia of Mathematics and Its Applications Vol. 132,
  Cambridge University Press, Cambridge (2011).

\bibitem{SchmidtWinter}
Schmidt G., Winter M.: Relational Topology.
LNM 2208 (2018).

\bibitem{Tarski}
Tarski A.: Introduction to Logic.
Oxford University Press (1941).

\bibitem{Winter}
Winter M.:
Strukturtheorie heterogener Relationenalgebren mit Anwendung auf Nichtdetermismus
in Programmiersprachen.
Dissertationsverlag NG Kopierladen GmbH, M\"unchen (1998)

\bibitem{Book}
  Winter M.:
  Goguen categories -- A categorical approach to $L$-fuzzy relations.
 Trends in Logic Vol. 25,
  Springer, Berlin (2007).

\bibitem{Fuzzy5}
Winter M.: Arrow Categories. Fuzzy Sets and Systems 160, 2893-2909 (2009).

\bibitem{Winter22-1}
Winter M.: Fixed Point Operators in Heyting Categories. Part I - Internal Fixed Point Theorem and Calculus.
(submitted to Journal of Pure and Applied Algebra, 2022).

\bibitem{PartI}
Winter M.: Relational Algebraic Approach to the Real Numbers - The Additive Group. Relational and Algebraic Methods in Computer Science, 20th International Conference, RAMiCS 2023,
LNCS 13896, 274-292 (2023)

\end{thebibliography}
\end{document}